\documentclass{ifacconf}

\usepackage{amsmath,amssymb,amsfonts,mathtools}
\usepackage{relsize}
\usepackage{natbib}            % you should have natbib.sty
\usepackage{graphicx}          % Include this line if your 
\usepackage[latex={-interaction=nonstopmode},crop=off,runs=2]{auto-pst-pdf} %use [off] to stop compilation
\usepackage{psfrag}
\newcommand{\D}{\mathbf {D}}
\newcommand{\Op}{\mathbf {Op}}
\newcommand{\Q}{\mathbf {Q}}
\newcommand{\U}{\mathbf {U}}
\newcommand{\X}{\mathbf {X}}
\newcommand{\Z}{\mathbf {Z}}
\newcommand{\T}{\mathcal {T}}
\newcommand{\B}{\mathcal {B}}
\newcommand{\A}{\mathcal {A}}
\newcommand{\V}{\mathcal {V}}
\newcommand{\G}{\mathcal {G}}

\newcommand{\la}{\langle}
\newcommand{\ra}{\rangle}

\newtheorem{definition}{Definition}
\newtheorem{example}{Example}

\begin{document}

\begin{frontmatter}

\title{Kernel representation approach to persistence of behavior} % Title, preferably not more than 10 words.

%\thanks[footnoteinfo]{This work was supported in part by the National Technological Agency. (sponsor and financial support acknowledgment goes here). Paper titles should be written in uppercase and lowercase letters, not all uppercase.}

\author[First]{Abdul Basit Memon} 
\author[Second]{Erik I. Verriest}

\address[First]{Georgia Institute of Technology, Atlanta, GA 30332 USA (abmemon@gatech.edu).}                                              
\address[Second]{Georgia Institute of Technology, Atlanta, GA 30332 USA (erik.verriest@gatech.edu)}

%\begin{keyword}                           % Five to ten keywords,  
%Cicero; Catiline; orations.               % chosen from the IFAC 
%\end{keyword}                             % keyword list or with the 
                                          % help of the Automatica 
                                          % keyword wizard

\begin{abstract}                          % Abstract of not more than 250 words.
The optimal control problem of connecting any two trajectories in a behavior $\B$ with maximal persistence of that behavior is put forth and a compact solution is obtained for a general class of behaviors. The behavior $\B$ is understood in the context of Willems's theory and its representation is given by the kernel of some operator. In general the solution to the problem will not lie in the same behavior and so a maximally persistent solution is defined as one that will be as close as possible to the behavior. A vast number of behaviors can be treated in this framework such as stationary solutions, limit cycles etc. The problem is linked to the ideas of controllability presented by Willems. It draws its roots from quasi-static transitions in thermodynamics and bears connections to morphing theory. The problem has practical applications in finite time thermodynamics, deployment of tensigrity structures and legged locomotion.
\end{abstract}
\end{frontmatter}

% Introduce the problem
% Why look into it
% Different areas where such connections can be found, tensigrity, thermo, gait transition
% Introduce previous work
% Highlight changes in work
% organization of paper

\section{Introduction}
The problem being considered here is that of connecting two trajectories from a set with a particular behavior in such a manner that the characteristic behavior persists during the transition. These particular behaviors could be stationary solutions, limit cycles or an even more general class of behaviors. The idea of exploring such transitions was first introduced in \cite{verriest2008first}. The problem can be stated as follows: Given two trajectories $w_1$ and $w_2$ of the same behavior, the objective is to construct a persistent transition, $w$, over some finite time interval $[a,b]$ such that $w=w_1$ for $t\leq a$ and $w=w_2$ for $t\geq b$. These ideas will be made more rigorous in the later sections. First we will motivate our interest in the problem of persistence of behavior and further elucidate the concept by some examples. The original motivation for the problem comes from the notion of quasi-static transitions in thermodynamics between two equilibrium points. Obviously if something is stationary then it cannot change but one can come arbitrarily close to the equilibrium points by slow motions. Thus persistence of stationarity is aimed for in this case. When such transitions are sought over finite time, it is closely related to our ideas  (\cite{berry2000book}, \cite{andresen1977step-carnot}). A related problem where such transitions are found is the deployment of tensigrity structures. In this case it is also desirable to transition from one configuration to another by remaining close to the equilibrium manifold, so that in case of loss of power the structure converges to some equilibrium configuration (\cite{sultan2003tensigrity}).

In the context of animal locomotion, gaits are periodic patterns of movement of the limbs. Most animals employ a variety of gaits such as one for walking and a different one for running (\cite{golubitsky2003symmetry}). To switch from a gait to another, one necessarily has to employ an aperiodic transition but animals do this naturally in a graceful manner. It is our hypothesis that this translates to the transient motion being as close as possible to periodic behavior. So the persistence of behavior required here will be periodicity. The theory of finding a persistent transition may also be of use in the control of legged robots (\cite{clark2006gaits}). A popular approach to legged robot control is to specify the gaits or different schemes of motion of a robot and then switch through these gaits. This reduces the complexity of the control problem. The problem then becomes one of finding a suitable gait transition that connects the two desired gaits from the set of dynamically consistent transitions. Thus the problem of finding a persistent transition is of significant practical interest. 

The problem of finding a persistent transition was presented in the earlier work: \cite{verriest2008first, deryck2009persistence, deryck2011thesis}.  %\footnote{Gluskabi is a mythical figure from the Penobscot Indians,who made himself from nothing, and transforms animals.} 
However the focus in the aforementioned papers was on specific behaviors. More general results are presented in this paper, which extend the earlier work in a number of ways. Firstly, a more generalized and rigorous mathematical formulation has been established and the nomenclature introduced in \cite{verriest2012mtns} is clarified. Secondly, the earlier Wronskian characterization of a scalar $n$-th order LTI differential system, introduced in \cite{verriest2012mtns}, is extended to the vector case. Thirdly, a very compact method is presented to find the transitions for a broad class of behaviors, characterized by the kernel of operators, with respect to any appropriate norm. This motivates the title of the paper. Fourthly, a similar compact representation for characterizing the transitions between trajectories of a linear time invariant dynamical system with respect to linear behaviors under any Sobolev norm has been found. Finally, all the ideas presented are illustrated using clear examples including one considering the optimal charging of a capacitor which is a significant problem in cyber-physics: the charging of batteries. 

The rest of the paper is organized as follows: A brief review of the behavioral approach of Willems is presented in Section \ref{sec:ba}. Following that in Section \ref{sec:gf}, the nomenclature and a mathematical formulation of our problem is presented. The Wronskian characterization in the vector case is presented in Section \ref{sec:wr}. The two main results of this paper are presented in Section \ref{sec:results} followed by some examples in Section \ref{sec:ex}. 

\section{Behavioral approach - A review}
\label{sec:ba}
We start by reviewing some of the relevant concepts from the behavioral approach to system theory. These ideas will be used later to set the nomenclature for our framework. A detailed exposition of the subject can be found in \cite{willems2007open,willemBook}. Let $\mathbb{T}$ denote the time axis. For continuous time systems we take $\mathbb{T}=\mathbb{R}$. %, and for discrete systems, $\mathbb{T}=\mathbb{Z}$. 
$\mathbb{W}$ is the set in which an $n$-dimensional observable signal vector, $w$, takes its values. Typically, $\mathbb{W}=\mathbb{R}^n$, $n\geq 1$. A dynamical system $\Sigma$ is defined as a triple
$\Sigma=(\mathbb{T},\mathbb{W},\mathcal{B})$. The behavior $\mathcal{B}$ is a suitable subset of $\mathbb{W}^\mathbb{T}$, for instance the piecewise smooth functions, compatible with the laws governing $\Sigma$. We define the evaluation functional $\sigma_t$ by $\sigma_t(w)=w(t)$ a.e. (exception where $w$ is not defined). The shift operator $\mathbf S_\tau$ is defined by $\sigma_t (\mathbf S_\tau w)= \sigma_{t+\tau}w$.

The dynamical system $\Sigma=(\mathbb{T},\mathbb{W},\mathcal{B})$ is said to be linear if $\mathbb{W}$ is a vector space over $\mathbb{R}$ or $\mathbb{C}$, and the behavior $\mathcal{B}$ is a linear subspace of $\mathbb{W}^\mathbb{T}$. The dynamical system $\Sigma=(\mathbb{T},\mathbb{W},\mathcal{B})$ is said to be shift invariant   if $w\in \mathcal{B}$ implies $\mathbf S_\tau w\in \mathcal{B}$
for all $\tau\in \mathbb{T}$. If $\Sigma=(\mathbb{T},\mathbb{W},\mathcal{B})$ is a shift-invariant dynamical system,  the behavior restricted to a small open interval $(-\epsilon,\epsilon)$ is defined by
$\mathcal{B}_\epsilon=\{\tilde{w}:(-\epsilon,\epsilon)\rightarrow \mathbb{W}\,|\,\exists w\in \mathcal{B}\; {\rm such}\;{\rm that}\;\sigma_t{\tilde{w}}=\sigma_t{w}\; {\rm for}\; {\rm all}\; -\epsilon<t<\epsilon\}.$ The continuous time system $\Sigma$ is called {\em locally specified} if for all $\epsilon>0$,\[(w\in \mathcal{B})\Leftrightarrow (\left.\mathbf S_\tau w\right|_{(-\epsilon,\epsilon)}\in \mathcal{B}_\epsilon\; {\rm for}\;{\rm all} \;\tau\in \mathbb{R}).\]
The behavior defined by the system of differential equations
\[R(\D)w=0,\quad R(\xi)\in \mathbb{R}^{p\times n}[\xi]\]
where $R(\xi)$ is a matrix of polynomials with real coefficients and $\D$ is the differentiation operator, represents a system of $p$ linear time invariant (LTI) ordinary differential equations (ODE) in $n$ scalar variables. A system described by behavioral differential equations is locally specified. In order to verify if a trajectory $w$ belongs to the behavior, it suffices to look at the trajectory
in an infinitesimal neighborhood about each point.

A behavior is called {\em autonomous} if for all $w_1,w_2\in\mathcal{B}$ $w_1(t)=w_2(t)  \; {\rm for}\; t\leq 0$ implies $w_1(t)=w_2(t)$ for almost all $ t$. For an autonomous system, the future is entirely determined by its past. The notion of {\em controllability} is an important concept in the behavioral theory. Let $\mathcal{B}$ be the behavior of a linear time invariant system.  This system is called controllable if for any two trajectories $w_1$ and $w_2$ in  $\mathcal{B}$, there exists a $\tau\geq 0$ and a trajectory $w\in \mathcal{B}$ such that
\[
\sigma_t(w)=\left\{
\begin{array}{ll}\sigma_t(w_1) & t\leq 0\\\sigma_t(\mathbf S_{-\tau}w_2)& t\geq \tau\end{array}\right.\]
i.e., one can switch from one trajectory to the other, with perhaps a delay, $\tau$. Note that an autonomous system cannot get off a trajectory once it is on it. Hence an autonomous system is not controllable. 

\section{Gluskabi framework}
\label{sec:gf}
In this section, we will first define the requisite nomenclature for our problem. We will then rigorously formulate our problem using the behavioral approach to system theory by Willems. We begin by defining a behavior which restricts the universum, $\mathbb W^\mathbb T$, which is the collection of all maps from the set of independent variables to the set of dependent variables, to just the ones which are interesting. \\

\begin{definition} The {\bf Base Behavior} ($\B_0$) is a subset of the universum $\B_0\subset \mathbb W^\mathbb T$ that defines the set of all allowable functions of interest. For any particular problem, the functions we are trying to connect lie in this set and the search for a connection\footnote{This usage of the term connection is different from a connection defined in differential geometry.} between the two is also conducted in this set. 
\label{def:base}
\end{definition}
For example, if we want to work with smooth functions entirely then $\B_0 = C^\infty(\mathbb T,\mathbb W)$. Or, if we are interested in the smooth trajectories of an LTI differential system then $\B_0 = \{w\in C^\infty(\mathbb T,\mathbb W) \;\text{s.t.}\; R(\D)w=0\}$, where $R(\xi)$ is a matrix of polynomials with real coefficients and $\D$ is the differentiation operator. In this paper we will be fixing our time axis and signal space to be real and so from this point onwards $\mathbb T = \mathbb R$ and $\mathbb W = \mathbb R^n$ for some $n\geq 1$. \\%The base system is our function universe and we will never step out of it. The trajectories we are trying to connect and the raccordations will all be in this behavior. 

\begin{definition} 
A \textbf{Type ($\T$)} is a strict subset of the base behavior ($\T\subset\B_0$) described by an operator $\Op:\A\to \mathcal V$ in the following way.
	\[ \mathcal T =\{w\in\mathcal A \;\text{such that}\; {\Op}\,w = 0 \}\]
where $\A\subset\B_0$ is the maximal linear space in the base behavior on which the operator is properly defined $\A\subset Dom(\Op)$ and $\mathcal V$ is a linear space as well. 
\label{def:type}
\end{definition}
The Type behavior defines the set of trajectories possessing a desired quality, which we want to connect. Given the obvious similarities, we call this the \emph{Kernel} representation of the type irrespective of whether the operator $\Op$ is linear or nonlinear. A type may admit representations other than the kernel representation but in this paper we will only consider the kernel representation of types. \\

\begin{definition}
A \textbf{Trait ($\T_\theta$)} is a subtype of the type i.e., it is a subset of the type such that it has its own characteristic behavior, given by some operator $\Op_\theta$. 
\[ \mathcal T_\theta =\{w\in\mathcal T \;\text{such that}\; {\Op}_{\theta}w = 0 \} \] 
\end{definition}
For instance, a trait could be specified by some (or all) boundary conditions, or some intermediate values and their derivatives. \\

\begin{example}[Constants]
Let $\B_0 = C^0(\mathbb R,\mathbb R)$ and $\D$ be the differential operator. Then the operator $\Op:=\D$ defines the type of constants in $\B_0$. An example of a particular trait in this type could be the constant $c$ i.e., $\T_c = \{w\in\T \;\text{s.t.}\; w = c\}$.
\end{example}

\begin{example}[Polynomials]
Let $\mathcal B_0 = C^0(\mathbb R,\mathbb R)$ and $\D$ be the differentiation operator. Then the operator $\Op:=\D^3$ defines the second order polynomials type in $\B_0$. An example of a trait in the second order polynomials type is the subtype of first order polynomials or constants. Another example of trait in this type is polynomials that vanish at $t=0$.
\end{example}

\begin{example}[Periodic signals with period $\tau$]
The operator $\Op:= (\mathbf I - \mathbf S_\tau)$ where $\mathbf I$ is the identity operator and $\mathbf S$ is the shift operator, defines the periodic type in $\B_0$. A smooth periodic function can be seen as a sum of harmonic signals of integer multiples of the base frequency. Thus the periodic type in $\B_0=C^\infty(\mathbb R,\mathbb R)$ may also be characterized by the infinite product operator $\left[\D\prod_{n=1}^\infty{\left(1+\frac{1}{n^2\omega^2}\D^2\right)}\right]$, which can also be written as $\sinh{\left(\frac{\pi}{\omega}\D\right)}$ (\cite{silverman1984ICA}), where $\omega=2\pi/\tau$. This representation defines a number of traits in terms of the number of finite product terms and these traits serve as various levels of approximation to the periodic functions. \\
\end{example}

The above three definitions form the basic nomenclature of our problem but we will need one more definition to rigorously define a connection later on. Given any type we can extend it to create a collection of related types in the following manner. \\

\begin{definition}
The \textbf{Equation Error System ($\T_{ee}$)} of a type $\T$, defined by the kernel of the operator $\Op$, is a union of behaviors $\T_e:= \left\{(w,e)\in\A\times \{e\} \;\text{s.t.}\; \Op\,w = e \right\}$.
\[ \mathcal T_{ee} := \cup_{e\in \mathcal V}\mathcal T_e = \left\{(w,e)\in\A \times \mathcal V \;\text{s.t.}\; \Op\,w = e \right\}\]
where $\mathcal V$ is the vector space where the image of $\Op$ lies i.e. $\Op(\A)\subset \mathcal V$. 
\end{definition}

Notice that the starting type $\T$ is the projection onto $\A$ of the behavior $\T_0 = \left\{(w,0)\in\A\times \V \;\text{s.t.}\; \Op\,w = 0 \right\}$ in this collection.  It is also worth noticing that the Equation Error System lies in an extended base behavior $\Sigma = (\mathbb T,\mathbb W\times\mathbb E,\B_0)$ where $\mathcal V\subset \mathbb E^\mathbb T$. 

\begin{example}
Consider the type in $C^\infty(\mathbb R,\mathbb R)$ defined by the operator $\Op:=(\D-\lambda \mathbf I)$ i.e., the type of multiples of the exponential $e^{\lambda t}$. Then the equation error system corresponding to this type is the set of solutions $w$ to the non-homogeneous ODE ($\D w-\lambda w = e$), for some forcing function $e\in C^\infty(\mathbb R,\mathbb R)$. \\
\end{example}

Now that we have suitable terminology we can formulate our problem. Given a type $\T$, the objective is to find a mapping which assigns to any two elements $w_1$ and $w_2$ in the said type a unique element $w$ in the base behavior which connects $w_1$ and $w_2$ in finite time i.e. over the interval $[a,b]$ and in such a manner that the defining quality of the type persists maximally. We will call this mapping the ``\emph{Gluskabi map}". Using the established idea that a type is given by the kernel of some operator $\Op$, the Gluskabi map and persistence of a trajectory is defined in the following manner. \\
\begin{definition}
Given a type $\T$ with the associated operator $\Op$, an element $w\in\A\subset\B_0$ is said to be \textbf{maximally persistent} with respect to the norm $\|.\|$, defined on $\mathcal V$ restricted to $[a,b]$, if $w$ minimizes $\|\Op\,w\|$.
\end{definition}
\begin{definition}
Given a type $\T$ with the associated operator $\Op$, the \textbf{Gluskabi map} $g:\T\times\T\to\B_0$ with respect to the norm $\|.\|$, defined on $\mathcal V$ restricted to $[a,b]$, is defined as follows 
\[g(w_1,w_2)(t)=\left\{\begin{array}{lr}w_1(t) & t\leq a \\ arg\!\min_{w\in\A}{\|\Op\,w\|} \quad & a<t<b\\ w_2(t) & t\geq b\end{array}\right.\]
\end{definition}

Clearly this requires that $\mathcal V$ restricted to the interval $[a,b]$ be a normed space. The connection in the interval $[a,b]$ will be called the ``\emph{Gluskabi raccordation}". As evident from the definition of the Gluskabi map, the element $w$ corresponding to $w_1,w_2\in\T$ may not lie in the type $\T$ and is constructed piecewise from elements in $\A$. A new behavior can now be constructed by collecting all the elements $w$ corresponding to any two elements $w_1$ and $w_2$ in the type $\T$ i.e., this behavior is the image of the Gluskabi map. This behavior will be called the ``\emph{Gluskabi Extension}" and can also be defined using the extended types $\T_{ee}$ in the following way. \\
\begin{definition}
Given a type $\T$ with the associated operator $\Op$, the \textbf{Gluskabi Extension ($\G_\T$)} with respect to the norm $\|.\|$, defined on $\mathcal V$ restricted to $[a,b]$, is defined as 
\begin{multline*} 
	\G_\mathcal T:= \left\{w\in\B_0 \;\text{s.t.}\; \exists w_1,w_2\in\T \text{ with }\Pi_{-}w = \Pi_{-}w_1,\right.\\\left. \Pi_{+}w = \Pi_{+}w_2,\text{and } \exists (u,e)\in\T_{ee} \right.\\ \left.\text{ s.t. } \Pi_{[a,b]}w = \Pi_{[a,b]}u \text{ with } \|e\|\text{ minimal}  \right\}
\end{multline*}
where $\Pi$ is the projection operator i.e. $\Pi_{-}w$ is the restriction of $w$ to the interval $(-\infty,a]$, $\Pi_{+}w$ is the restriction to the interval $[b,\infty)$ and $\Pi_{[a,b]}w$ is the restriction to the interval $[a,b]$. 
\end{definition}

Notice that $\T\subset\G_\T$. Since the space $\V$ generally admits multiple norms, the Gluskabi map and extension will in general depend on the chosen norm and the raccordation interval. Thus, a suitable norm in conjunction with the operator $\Op$ completely characterizes the desired persistence. For instance, if $\Op$ is a differential operator of some order acting on functions then any Sobolev norm of compatible degree can be used to get the required level of smoothness. Say the time interval is $[a,b]$ and the $\Op:C^r(\mathbb R,\mathbb R)\to C^s(\mathbb R,\mathbb R)$, then the Sobolev norm $\|.\|_W$ on $e\in \V=C^s([a,b],\mathbb R)$ is given by
\begin{multline*}
	\|e\|_W = \sum_{i=0}^n{\rho_i \|\D^i e\|_{L^2}} \quad \text{where }\rho_i> 0, \; n\leq s \\ \text{and } \|x\|_{L^2}^2 = \int_a^b{x^2(t) dt} 
\end{multline*}

\section{LTID Type}
\label{sec:wr}
In this section, we focus our attention on an interesting type namely the linear time invariant differential (LTID) behavior $\mathcal{L}^k_n$ of some order $n$ , i.e., the set of all solutions to any system of $k$ constant coefficient homogeneous differential equations of $n$th order. The goal here is to find a kernel representation for this type $\mathcal{L}^k_n$ i.e., find the operator characterizing this type in accordance with Definition-\ref{def:type}. This type was first introduced in \cite{verriest2012mtns}, where the operator was derived for the scalar $n$-th order differential equation case i.e., when $k=1$. Using Willems's approach, this behavior is represented as,
\begin{multline*}
	\mathcal{L}^k_n = \left\{ w\in C^n(\mathbb R,\mathbb{R}^{k})\left|\right. \exists R\in\mathbb{R}[\xi]^{k\times k}\right.\\	\left.\text{for which}\; R(\D)w=0\right\} 
\end{multline*}
where $\D$ is the differentiation operator and $R$ is a polynomial matrix
\[R(\xi):= R_0 \xi^n + R_1\xi^{n-1} + \cdots + R_n \xi .\]
Let's assume that $R_0=I$ and that the system of differential equations is not underdetermined or overdetermined. If $w\in\mathcal{L}^k_n$ then there exist $R_i\in\mathbb{R}[\xi]^{k\times k}$ such that the following holds true 
\begin{equation} w^{(n)} + R_1w^{(n-1)} + \cdots + R_n w = 0 \end{equation}
\begin{equation} \Rightarrow \left(\D^n+R_1\D^{n-1}\cdots+R_n\right) \left[\begin{array}{cccc} w & \dot{w} & \cdots & w^{(nk+k-1)} \end{array}\right] = 0 \end{equation}
\begin{equation}\setlength{\arraycolsep}{3pt} \begin{bmatrix}R_n &\cdots& R_1 & I\end{bmatrix}
	\begin{bmatrix}
		w & \dot{w} & \cdots & w^{(nk+k-1)} \\
		\vdots  & \vdots		& \ddots & \vdots	\\
		w^{(n)}       & w^{(n+1)}		& \cdots & w^{(n+nk+k-1)}
	\end{bmatrix} = \mathlarger{0} \label{eq:ltid}\end{equation}
Notice that the matrix on the right looks like a Wronskian in the vector functions $(w,\dot{w},\cdots,w^{(nk+k-1)})$. But if the vector $w$'s are expanded then this will not be a symmetric matrix. Let's call it the generalized Wronskian nonetheless and partition it in the following manner.
\begin{equation}
	 \arraycolsep = 2pt	 
	 \left[\begin{array}{ccc|ccc}
			w & \cdots & w^{(nk-1)} & w^{(nk)}& \cdots & w^{(nk+k-1)}\\
			\vdots &  & \vdots & \vdots &  & \vdots \\
			w^{(n-1)} & \cdots & w^{(n-1+nk-1)} & w^{(n-1+nk)} & \cdots & w^{(n+nk+k-2)}\\
			\hline
			w^n & \cdots & w^{(n-1+nk)} & w^{(n+nk)} & \cdots & w^{(n+nk+k-1)}
		\end{array}\right] \label{eq:partition}\end{equation}				
Let's name the upper left and the upper right block of this partitioned matrix as $\widehat{W}$ and $\widetilde{W}$ respectively. i.e., 
\begin{eqnarray} \widehat{W} &=& \begin{bmatrix} w & \cdots & w^{(nk-1)} \\ \vdots&&\vdots \\ w^{(n-1)} & \cdots & w^{(n-1+nk-1)}\end{bmatrix} \label{eq:name1}\\
\widetilde{W} &=& \begin{bmatrix} w^{(nk)}& \cdots & w^{(nk+k-1)} \\ \vdots&&\vdots \\ w^{(n-1+nk)} & \cdots & w^{(n+nk+k-2)}\end{bmatrix}. \label{eq:name2}\end{eqnarray}

Note that:
\[ \begin{bmatrix} I&O \\ -BA^{-1}&I \end{bmatrix}\begin{bmatrix} A&C \\ B&D \end{bmatrix}\begin{bmatrix} I&-A^{-1}C \\ O&I \end{bmatrix} = \begin{bmatrix} A&O \\ O&Schur(A) \end{bmatrix} \]
where $Schur(A)$ is the Schur complement of $A$. Using this fact, (\ref{eq:ltid}) can be written as 
\begin{multline} \arraycolsep = 2pt	 
	 \begin{bmatrix}R_n &\cdots& R_1 & I\end{bmatrix}
	 \left[\begin{array}{c|c} I & O \\ \hline \begin{bmatrix}w^n & \cdots & w^{(n-1+nk)}\end{bmatrix}\widehat{W}^{-1}& I \end{array}\right] \\
	 \left[\begin{array}{c|c} \hat{W} & O \\\hline O & Schur(\widehat{W}) \end{array}\right] = \mathlarger{0} \end{multline}

\begin{equation}\Rightarrow\left\{\begin{array}{l} 
	\begin{bmatrix}R_n &\cdots& R_1\end{bmatrix}\widehat{W}+\begin{bmatrix}w^n & \cdots & w^{(n-1+nk)}\end{bmatrix} = 0 \\
	Schur(\widehat{W}) = 0 
\end{array} \right. \label{eq:lti_final}\end{equation}

The first equation in (\ref{eq:lti_final}) is just a subset of the original set of equations (\ref{eq:ltid}), specifically the ones formed by using the columns to the left of the partition in (\ref{eq:partition}). Thus, if $w\in\mathcal{L}^k_n$ then a necessary condition for $w$ is that $Schur(\widehat{W})=0$ i.e.,
\begin{multline}
	 \begin{bmatrix}w^{(n+nk)} & \cdots & w^{(n+nk+k-1)} \end{bmatrix} - \\
	 \begin{bmatrix}w^n & \cdots & w^{(n-1+nk)}\end{bmatrix} 
	 \widehat{W}^{-1}\widetilde{W} = \mathlarger{0}
	\label{eq:ltiop}
\end{multline}
where $\widehat{W}$ and $\widetilde{W}$ are as defined in (\ref{eq:name1}) and (\ref{eq:name2}). Thus we have found a nonlinear operator $\Op$ such that the functions $w$ satisfying $\Op\,w=0$ or (\ref{eq:ltiop}) form the $n$th order LTID type in $k$ variables ($\mathcal L^k_n$). 

\section{Finding the Gluskabi extension}
\label{sec:results}
Now that we have rigorously stated our problem, we will devote this section to present two results on finding the Gluskabi extension that are applicable to a broad collection of types. The only requirement for these results to be applicable is that the range of the operator associated with the type or $\V$, restricted to the raccordation interval $[a,b]$ be an inner product space and the said operator admits an adjoint. This condition is not extremely restrictive and is satisfied by a number of interesting operators such as differential operators and shift operators. The two results differ in the choice of the base behavior; in the first result the base behavior is some appropriately chosen function space where as in the second case the base behavior is the set of smooth trajectories of a dynamical system. So the two cases are appropriately called the \emph{signal raccordation} and the \emph{dynamical raccordation} problem respectively. 

\subsection{Signal Raccordation}
\begin{thm}
Given a type $\T$ with the associated operator $\Op$, the Gluskabi extension with respect to the norm $\|.\|_\Q$, where $\Q$ is a self-adjoint operator, is given by
\[ \G_\T|_{[a,b]} = \left\{w\in\B_0 \;\text{ such that }\; \Op_w^*\Q\Op\,w = 0 \right\}|_{[a,b]},\]
where the raccordation is sought over the interval $[a,b]$, the norm is computed as $\|.\|^2_\Q=\la\Q(.),(.)\ra$ and $\Op_w$ is the linearized form (G\^{a}teaux derivative) of the operator $\Op$ about $w$.
\label{th:sig}
\end{thm}

\begin{pf}
This can be easily proved using variational calculus. Given the type operator $\Op:\A\to\V$ and the norm $\|.\|^2_\Q$, the cost functional to be minimized can be written as
\begin{equation} J(w) = \|\Op w\|^2_\Q = \la \Q\Op\, w, \Op\, w\ra \label{eq:cost_s}\end{equation}
Now using the assumption that $\Op$ is G\^{a}teaux differentiable, it is shown that the first variation of $J$ exists and its expression in terms of $\Op$ is computed as follows:
\begin{align}
\Delta J=J(w+th)-J(w) &= \la\Q\Op (w+th),\Op (w+th)\ra \nonumber\\
&\quad - \la\Q\Op\, w,\Op\, w\ra \nonumber
\end{align}
%We know that 
It is given that
$\Op (w+th) = \Op\, w+\Op_w\, th+O\left(t^2\right)$. So,
\begin{multline}
\Delta J= 2t\la\Q\Op\, w,\Op_w\, h\ra + t^2\la\Q\Op_w\, h,\Op_w\, h\ra \\
+ 2\la\Q\Op\, w,O\left(t^2\right)\ra + 2\la\Op_w\, th,O\left(t^2\right)\ra \\
+ \la O\left(t^2\right),O\left(t^2\right)\ra. 
\label{eq:delJ}
\end{multline}
The last expression is obtained using the facts that $\Q$ is self adjoint and $\Op_w$ is linear. Then the first variation is given by
\begin{align}
\delta J(w;h) &= \lim_{t\to 0}{\frac{\Delta J}{t}} \nonumber \\
&= 2\la\Q\Op\, w,\Op_w\, h\ra \nonumber \\
&= 2\la\Op_w^*\Q\Op\, w, h\ra + \text{boundary terms}
\end{align}
since each of the other terms in (\ref{eq:delJ}) goes to zero as $t\to 0$. If $w$ is the minimizer of the functional (\ref{eq:cost_s}) then the first variation $\delta J(w;h)$ is zero at $w$. Thus a necessary condition for all raccordations in the Gluskabi extension $\G_\T$ is that 
\begin{equation}
	\Op_w^*\Q\Op\, w = 0 \quad\forall w\in \G_\T.
	\label{eq:glusk_sig}
\end{equation}
The boundary terms are zero because of the given boundary conditions for the problem i.e., $w$ and possibly a number of its derivatives at $t=a$ and $t=b$ are fixed. Therefore the admissible variations $h$ are zero at the endpoints. \qed
\end{pf}

If there exists an operator $\Op^*$ such that 
\begin{equation}
	\Op^*(w+ \delta w)-\Op^*\,w = \Op_w^*\,\delta w\quad \forall w\in\A 
\end{equation}
then the above condition for the Gluskabi Extension (\ref{eq:glusk_sig}) can be written as the following nested form.
\begin{equation}
	\Op^*(w+ \Q\Op\,w) = \Op^*\,w \quad\forall w\in \G_\T
\end{equation}
Furthermore an example of the norms that can be employed is the Sobolev norm. For this case the operator corresponding to the norm is $\Q=W(-\D^2)$, where $W(\xi)=1+\xi+\xi^2+\cdots+\xi^k$ with $k$ being order of the norm and $\D$ is the differentiation operator. 

\subsection{Dynamical Raccordation}
Next we look at the dynamical raccordation case when the trajectories in the base behavior are constrained by the dynamics of the system. Since one is never allowed to step out of the base behavior, we can call the dynamical system constraints ``\emph{hard constraints}" where as the type constraints are ``\emph{soft constraints}". The focus of the following result is on finding the Gluskabi extension for polynomial differential types i.e. $\Op$ is a polynomial in $\D$ and the base behavior is trajectories of an LTI dynamical system i.e. $\B_0=\{w\in C^\infty(\mathbb R,\mathbb R^q)\;\text{s.t.}\;R(\D)\,w=0\}$ where $R\in\mathbb R^{g\times q}[\xi]$ and $g<q$. 

The presentation of the main result is preceded by some necessary remarks. Given a scalar type $(\T,\Op)$ i.e., defined on signal space $\mathbb W=\mathbb R$, it can be correspondingly defined for vector trajectories i.e., $\mathbb W=\mathbb R^q$ by extending $\Op$ as $\Op^e\,w=\left(\Op\,w_1,\cdots,\Op\,w_q\right)^T$ for $w=(w_1,\cdots,w_q)^T\in\mathbb W^\mathbb R$. In the following result $\Op$ will be understood to be $\Op^e$ wherever appropriate. The inner product is appropriately extended as well. Every LTI system has an equivalent minimal representation that can also be expressed in the input/output form $P(\D)y = N(\D)u$ where $P\in\mathbb R^{g\times g}[\xi]$, $det\,P\neq 0$, and $P^{-1}N$ is a proper matrix (\cite{willemBook}). This input/output form of an LTI system will be used in the following result and since $u$ and $y$ are simply obtained by some partition of $w$, $w\in\T$ implies that both $u$ and $y$ are of the same type. Hence we are looking for connections of input/output pairs of the type $\T$. \\

\begin{thm}
Given a minimal and controllable linear time invariant dynamical system $P(\D)y = N(\D)u$ and a type $\T$ with the associated linear operator $\Op$, the trajectories in the Gluskabi extension with respect to the Sobolev norm $\|.\|_\Q$, restricted to the interval $[a,b]$ is given by the following equations:
\begin{align*}
\left(\U_{12}^*\,{\Op^u}^*\,\Q^u\,\Op^u\, \U_{12} + \U_{22}^*\,{\Op^y}^*\,\Q^y\,\Op^y\, \U_{22}\right)\,\eta &= 0 \\
-\U_{12}\,\eta &= u \\
\U_{22}\,\eta &= y 
\end{align*}
where $U = \begin{bmatrix}U_{11} & U_{12} \\ U_{21} & U_{22}\end{bmatrix}$ is a unimodular matrix such that $\begin{bmatrix}N & P\end{bmatrix} U = \begin{bmatrix} I & O \end{bmatrix}$ and $P\in\mathbb R^{g\times g}[\xi]$, $N\in\mathbb R^{g\times (q-g)}[\xi]$, $U\in\mathbb R^{q\times q}[\xi]$, $\Q$ is self-adjoint, $y$ and $u$ are the output and input respectively and $\D$ is the differentiation operator. 
\label{th:lti-cont}
\end{thm}

\begin{pf}
The cost function to be minimized along with the adjoined constraints is given in the inner product form:
\begin{multline}
J(u) = \frac{1}{2}\la \Q^u\Op^u\, u,\Op^u\, u\ra + \frac{1}{2}\la\Q^y\Op^y\, y,\Op^y\, y\ra \\
+ \la\lambda,P(\D)y - N(\D)u\ra 
\end{multline}
where $\Op^u$, $\Op^y$, $\Q^u$ and $\Q^y$ are the appropriately extended forms of the operators $\Op$ and $\Q$ depending on dimensions of $u$ and $y$ respectively. The first variation of the cost function due to a perturbation in $u$ can be computed as follows,
\begin{multline}
 \delta J(u;\delta u) = \la{\Op^u}^*\Q^u\Op^u\, u,\delta u\ra + \la{\Op^y}^*\Q^y\Op^y\, y,\delta y\ra \\
+ \la P(\D)^*\,\lambda,\delta y\ra  - \la N(\D)^*\,\lambda,\delta u\ra \\
+ \text{ boundary terms}
\end{multline}
where we have used the fact that $\Q^u$ and $\Q^y$ are self adjoint, since $\Q=\Q^*$, in deriving the last expression. The boundary terms can be ignored since the functions $u$ and $y$ over the interval $[a,b]$ are to be matched to their respective given trajectories at the boundaries. Thus the variations $\delta u$ and $\delta y$ and appropriate number of their derivatives are zero at the end points. This leads to the Euler-Lagrange equations, 
\begin{equation}
	{\Op^y}^*\Q^y\Op^y\, y + P(\D)^*\,\lambda = 0.
	\label{eqn:el}
\end{equation}
The necessary condition for optimality is,
\begin{equation}
	{\Op^u}^*\Q^u\Op^u\, u - N(\D)^*\,\lambda = 0
	\label{eqn:oc}
\end{equation}
To find the Gluskabi extension it is required to eliminate $\lambda$ from the above two equations and solve the resultant equations along with the dynamical system equation for $u$ and $y$ and the given boundary conditions. In other words one needs to find the behavior given by the representation,
\begin{equation}
\begin{bmatrix} N^* & X & O \\ P^* & O & Z \\ O & N & P \end{bmatrix}(\D)
\begin{bmatrix} \lambda \\ -u \\ y \end{bmatrix} 
= 0
\label{eq:ham_sys}
\end{equation}
where $X$,$Z$,$N^*$, and $P^*$ are polynomial matrices such that $X(\D)={\Op^u}^*\Q^u\Op^u$, $Z(\D)={\Op^y}^*\Q^y\Op^y$, $N^*(\D)=N(\D)^*$, and $P^*(\D)=P(\D)^*$. This behavior in (\ref{eq:ham_sys}) will be unchanged under any left unimodular transformation on the polynomial matrix (\cite{willemBook}). Since the system is controllable, the rank of the matrix $\begin{bmatrix}P(s) & -N(s)\end{bmatrix}$ is the same for all $s\in\mathbb C$ and because of minimality the matrix has full row rank for almost all $s$. This implies that this matrix has full rank for all $s$ and the polynomial matrices $P$ and $N$ are left coprime (\cite{willemBook} and \cite{kailathbook}). Thus there  always exists a unimodular matrix $U$ such that 
\begin{equation}
	\begin{bmatrix} N & P \end{bmatrix} U = \begin{bmatrix} I & O \end{bmatrix}. 
\end{equation}
It also holds that
\begin{equation}
	U^* \begin{bmatrix} N^* \\ P^*\end{bmatrix}  = \begin{bmatrix} I \\ O \end{bmatrix}
\end{equation}
where $U^*(s)=U(-s)^T$. The matrix $U^*$ is also unimodular since $det\,U = det\,U^T$ and since the determinant is a polynomial in the entries of the matrix, which are themselves polynomials for the matrix $U$ and so the determinant is a polynomial in the indeterminate $s$ and is some constant since $U$ is unimodular and so changing the indeterminate to $-si$ doesn't change the determinant. A new unimodular matrix can now be constructed using $U^*$ which is $\begin{bmatrix}U^* & O \\ O & I\end{bmatrix}$ and applying it to polynomial matrix in (\ref{eq:ham_sys}) yields
\begin{equation}
\begin{bmatrix}U^* & O \\ O & I\end{bmatrix}
\begin{bmatrix}N^* & X & O \\ P^* & O & Z \\ O & N & P \end{bmatrix}.
\end{equation}
If the matrix $U$ is partitioned as $\begin{bmatrix}U_{11} & U_{12} \\ U_{21} & U_{22}\end{bmatrix}$ then $U^* = \begin{bmatrix}U_{11}^* & U_{21}^* \\ U_{12}^* & U_{22}^*\end{bmatrix}$ and the above expression simplifies to:
\begin{equation}
\begin{bmatrix}I & U_{11}^* X & U_{21}^* Z \\ O & U_{12}^* X & U_{22}^* Z \\ O & N & P \end{bmatrix}.
\end{equation}
The behavior corresponding to the polynomial matrix above is equivalent to the one in (\ref{eq:ham_sys}).
\begin{align}
\begin{bmatrix}I & U_{11}^* X & U_{21}^* Z \\ O & U_{12}^* X & U_{22}^* Z \\ O & N & P \end{bmatrix}
\begin{bmatrix}I & O \\ O & U \end{bmatrix}
\begin{bmatrix}I & O \\ O & U \end{bmatrix}^{-1}(\D)
\begin{bmatrix} \lambda \\ -u \\ y \end{bmatrix} &= 0 \nonumber\\
\begin{bmatrix}I & \U_{11}^* \X \U_{11}+\U_{21}^* \Z \U_{21} & \U_{11}^* \X \U_{12}+ \U_{21}^* \Z \U_{22}\\ O & \U_{12}^* \X \U_{11}+ \U_{22}^* \Z \U_{21} & \U_{12}^* \X \U_{12} + \U_{22}^* \Z \U_{22}\\ O & I & O \end{bmatrix}\!\!
\begin{bmatrix} \lambda \\ \nu \\ \eta \end{bmatrix} \!&=\! 0 
\label{eq:uni}
\end{align}
where the bold font corresponds to the differential operator of the respective polynomial e.g. $\X=X(\D)$ and so on, $\begin{bmatrix}I & O \\ O & U \end{bmatrix}^{-1}\begin{bmatrix} \lambda \\ -u \\ y \end{bmatrix} = \begin{bmatrix} \lambda \\ \nu \\ \eta \end{bmatrix}$, $\nu$ is a $g\times 1$ vector, and $\eta$ is a $(q-g)\times 1$ vector. The third row of (\ref{eq:uni}) simplifies to $\nu = 0$ and the second row simplifies to the equation,
\[ \left(\U_{12}^*\X \U_{12} + \U_{22}^*\Z \U_{22}\right)\,\eta = 0 \]
or
\begin{equation}
	\left(\U_{12}^*\,{\Op^u}^*\,\Q^u\,\Op^u\, \U_{12} + \U_{22}^*\,{\Op^y}^*\,\Q^y\,\Op^y\, \U_{22}\right)\,\eta = 0 
\end{equation}
and the substitution yields
\begin{align}
u &= -\U_{12}\,\eta \\
y &= \U_{22}\,\eta 
\end{align}
\begin{flushright}
\qed
\end{flushright}
\end{pf}

It is of course assumed that the set $\{(u,y) \text{ s.t. } P(\D)y = N(\D)u \text{ and } \Op\,u=0 \text{ and } \Op\,y=0\}$ is nonempty i.e., there exist input/output pairs of the dynamical system of the required type $\T$. Otherwise the question of finding the Gluskabi extension is moot. Also the controllability assumption (as defined in Section \ref{sec:ba}) is a sufficient condition for the solution to exist. It guarantees that there exist trajectories of the dynamical system connecting the left trajectory to the right one in some finite time. Furthermore, for smooth solutions to an LTID system the time can be taken to be arbitrarily small (\cite{willemBook}) and so the length of the interval $[a,b]$ does not matter. This result can be further generalized to the case when only the input or the output is of the type and needs to be connected or to the case when the persistence of output is more important then the input. Either of these cases can be viewed as an extension of the previous result by changing the operator $\Q$ of the inner product. For instance, in the first case can be accomplished by choosing $\Q^u=0$ or $\Q^y=0$. This will be further elucidated in the final example in Section \ref{sec:ex}.

\section{Examples}
\label{sec:ex}
\begin{figure}[h]
\centering
\psfrag{wl}{$5e^{-2t}$}
\psfrag{wr}{$0.02e^{8t}$}
\includegraphics[scale=0.5]{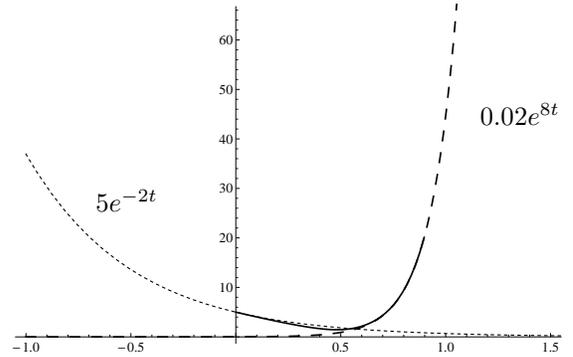}
\caption{A raccordation (solid line) connecting $w_1 = 5e^{-2t}$ (dotted line) in the $\mathcal L^1_1$ type to $w_2 = 0.02e^{8t}$ (dashed line) in the same type.}
\label{fig:lti}
\end{figure}

In this section, we illustrate the results presented in the previous section with the help of some examples. We start by looking at the signal raccordation problem for which the result was presented in Theorem \ref{th:sig}. Let's choose our base behavior to be $\B_0=C^0(\mathbb R,\mathbb R)$ and the type to be scalar first order LTID type $\mathcal L^1_1$ i.e. the set of all exponentials $ce^{\lambda t}$ for all values of $c\in\mathbb R$ and $\lambda\in\mathbb R$. Looking back at Section \ref{sec:wr}, the operator for this type is found to be $\Op\,w = \ddot{w}w - \dot{w}^2$. Say the raccordations are sought over the interval $[0,1]$ and the norm to be minimized is the usual $L^2$ norm. Then according to Theorem \ref{th:sig}, the raccordation $w$ over the interval $[0,1]$ must be the solution to the differential equation $\Op_w^*\Op\,w = 0$ where $\Op_w = w\D^2-2\dot{w}\D+\ddot{w}I$ or 
\[ w^{(4)}w^2+2w^{(3)}\dot{w}w-3\ddot{w}\dot{w}^2 = 0 \]
This gives us a generalized solution and then the specific raccordation connecting say $w_1$ and $w_2$ is obtained by using the boundary conditions i.e. 
$w^{(i)}(0)=w_1^{(i)}(0)$ and $w^{(i)}(1)=w_2^{(i)}(1)$ for $i=0$ and $i=1$. The raccordation for the case when $w_1 = 5e^{-2t}$ and $w_2 = 0.02e^{8t}$ is shown in Fig. \ref{fig:lti}.

\begin{figure}[h]
\centering
\psfrag{u}{$u(t)$}
\psfrag{y}{$y(t)$}
\includegraphics[scale=0.5]{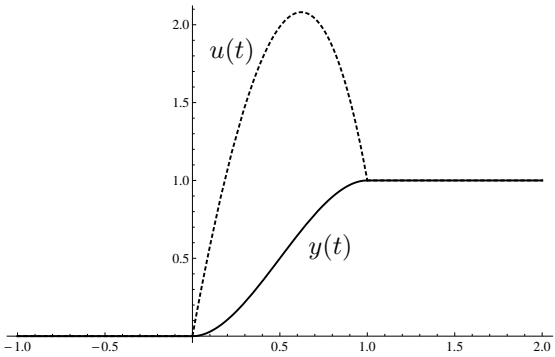}
\caption{The raccordation from constant $0$ to the constant $1$. The input is dashed line and output is the solid one.}
\label{fig:dy}
\end{figure}

Next we look at an example for the dynamical raccordation case of Theorem \ref{th:lti-cont}. We have a scalar first order LTI system given by the input-output differential equation $(\D+1)y = u$. We are interested in transitioning from one constant steady state to another. So our type is constants and $\Op=\D$. Notice that elements of this type satisfy the hard constraint i.e., if $y = c$ where $c$ is some constant then $u = c$. The transfer function for this system is $H(s) = \frac{1}{s+1}$ and so at steady state $y_{ss} = u_{ss}$, by the final value theorem. The chosen norm is again the $L^2$ norm and the raccordation time interval is $[0,1]$. The numerator and denominator polynomials are $N(s)=1$ and $P(s)=s+1$ respectively. And so $U = \begin{bmatrix}1 & -(s+1) \\ 0 & 1 \end{bmatrix}$ is the unimodular matrix required by Theorem \ref{th:lti-cont} and to find the Gluskabi extension the following system of equations need to be solved. 
\begin{align}
[(\D+1)^*\Op^*\Op(\D+1)+\Op^*\Op]\,\eta &= 0 \label{eq:ex_eta}\\
(\D+1)\,\eta &= u \\
\eta &= y
\end{align}
The equation (\ref{eq:ex_eta}) is simplified to get
\begin{equation}
	\left(\D^4-2\D^2\right)\,\eta = 0 
\end{equation}
Solving these differential equations yields,
\begin{align}
y(t) &= A e^{\sqrt{2}t} + B e^{-\sqrt{2}t} + C + Dt \\
u(t) &= (1+\sqrt{2})A e^{\sqrt{2}t} + (1-\sqrt{2})B e^{-\sqrt{2}t} + C + D(1+t)
\end{align}
Again the specific raccordation is obtained by using the boundary conditions i.e. $u(0)$, $y(0)$, $u(1)$, and $y(1)$. The raccordation for the case when $u=y = 0 \text{ for } t\leq 0$ and $u=y = 1 \text{ for } t\geq 1$ is illustrated in Fig. \ref{fig:dy}. \\ 

\begin{figure}[h]
\centering
\psfrag{u}{$u(t)$}
\psfrag{R}{$R$}
\psfrag{C}{$C$}
\includegraphics[scale=0.5]{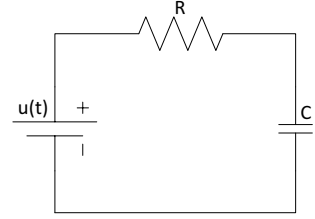}
\caption{Charging of a capacitor in an RC circuit}
\label{fig:circuit}
\end{figure}

\begin{figure}[h]
\centering
\psfrag{u}{\small $u(t)$}
\psfrag{c}{\small $q(t)$}
\psfrag{QD}{\footnotesize $Q$}
\psfrag{QC}{\footnotesize $\frac{Q}{C}$}
\psfrag{T}{\scriptsize $T$}
\psfrag{t}{\scriptsize $t$}
\psfrag{q}{\scriptsize $q,u$}
\includegraphics[scale=0.5]{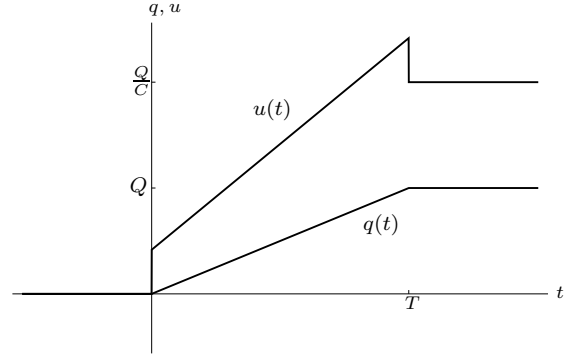}
\caption{Charge and input voltage trajectories of the RC circuit}
\label{fig:capacitor}
\end{figure}

We end this section by looking at the cyber-physical problem of charging a capacitor. We consider the simplest series RC circuit shown in Fig. \ref{fig:circuit}. The objective here is to put a charge $Q$ on the capacitor in time interval $[0,T]$. So the type to be considered for this case is the type of constants and again the $L^2$ norm is minimized. The dynamical system equation associated with the circuit is $\dot{q}+\frac{1}{RC}q = \frac{1}{R}u$, where $q$ is the charge on the capacitor and $u$ is the source voltage as well as the input over here. The type constraint is only imposed on the output i.e. $q$ and so in terms of Theorem \ref{th:lti-cont}, $\Q^u=0$. The resulting trajectory of charge and the input voltage is illustrated in Fig. \ref{fig:capacitor}. An interesting parallel has been found that the resulting minimizing trajectory obtained from applying Theorem \ref{th:lti-cont} is the same trajectory obtained when minimizing the heat generated in the resistor as shown in \cite{vos2000equipartition}. This points to a possible correlation between our theory and minimization of entropy for thermodynamic systems and will be explored in future publications. 

\section{Conclusion}
The previous work of Verriest and Yeung was extended by introducing new terminology and rigorously formulating the raccordation problem using those terms. The solution to the raccordation problem corresponds to constructing the Gluskabi Extension. A generalized construction of the Gluskabi Extension was obtained for the class of types defined by the kernel of some operator, which admits an adjoint. The Gluskabi Extension for linear types constrained by the trajectories of a linear time invariant dynamical system was also obtained. Finally a novel operator characterization for the LTI $n$th order differential type was developed as well.

\bibliography{gluskabi}           % and a bib file to produce the 

                                 % bibliography (preferred). The
                                 % correct style is generated by
                                 % Elsevier at the time of printing.

\end{document}